\documentclass{amsart}
\usepackage{amsopn}
\usepackage{amssymb, amscd}

\newcommand{\nc}{\newcommand}

\nc{\vg}{\mathfrak{v} } \nc{\wg}{\mathfrak{w} }
\nc{\zg}{\mathfrak{z} } \nc{\ngo}{\mathfrak{n} }
\nc{\ngoq}{\mathfrak{n}^{\QQ} } \nc{\ngoz}{\mathfrak{n}^{\ZZ} }
\nc{\ggoq}{\mathfrak{g}^{\QQ} } \nc{\kg}{\mathfrak{k} }
\nc{\mg}{\mathfrak{m} } \nc{\bg}{\mathfrak{b} }
\nc{\ggo}{\mathfrak{g} } \nc{\ggob}{\overline{\mathfrak{g}} }
\nc{\sog}{\mathfrak{so} } \nc{\sug}{\mathfrak{su} }
\nc{\spg}{\mathfrak{sp}} \nc{\slg}{\mathfrak{sl} }
\nc{\glg}{\mathfrak{gl} } \nc{\cg}{\mathfrak{c} }
\nc{\rg}{\mathfrak{r} } \nc{\hg}{\mathfrak{h} }
\nc{\tg}{\mathfrak{t} } \nc{\ug}{\mathfrak{u} }
\nc{\dg}{\mathfrak{d} } \nc{\ag}{\mathfrak{a} }
\nc{\pg}{\mathfrak{p} } \nc{\sg}{\mathfrak{s} }
\nc{\lgo}{\mathfrak{l} } \nc{\fg}{\mathfrak{f} }

\nc{\pca}{\mathcal{P}} \nc{\nca}{\mathcal{N}}
\nc{\lca}{\mathcal{L}} \nc{\oca}{\mathcal{O}}
\nc{\mca}{\mathcal{M}} \nc{\tca}{\mathcal{T}}
\nc{\aca}{\mathcal{A}}

\nc{\vp}{\varphi} \nc{\ddt}{\frac{{\rm d}}{{\rm d}t}}
\nc{\im}{\mathtt{i}}

\renewcommand{\l}{\lambda}
\nc{\ala}{Anosov Lie algebra} \nc{\alas}{Anosov Lie algebras}

\nc{\SO}{{\mathrm SO}} \nc{\Spe}{{\mathrm Sp}} \nc{\Sl}{{\mathrm
SL}} \nc{\SU}{{\mathrm SU}} \nc{\Or}{{\mathrm O}} \nc{\U}{{\mathrm
U}} \nc{\Gl}{{\mathrm GL}} \nc{\Se}{{\mathrm S}} \nc{\Cl}{{\mathrm
Cl}} \nc{\Spin}{{\mathrm Spin}} \nc{\Pin}{{\mathrm Pin}}

\nc{\RR}{{\mathbb{R}}} \nc{\HH}{{ \mathbb{H}}} \nc{\CC}{{
\mathbb{C}}} \nc{\ZZ}{{ \mathbb{Z}}} \nc{\FF}{{ \mathbb{F}}}
\nc{\NN}{{ \mathbb{N}}} \nc{\QQ}{{ \mathbb{Q}}} \nc{\PP}{{
\mathbb{P}}}

\nc{\vs}{\vspace{.5cm}}
 \nc{\ip}{\langle\cdot,\cdot\rangle}
\nc{\la}{\langle} \nc{\ra}{\rangle} \nc{\unm}{\frac{1}{2}}
\nc{\unc}{\frac{1}{4}} \nc{\und}{\frac{1}{16}} \nc{\f}{\frac}

\nc{\no}{\vs\noindent} \nc{\lam}{\Lambda^2\ggo^*\otimes\ggo}
\nc{\tang}{{\rm T}} \nc{\dif}{{\rm d}} \nc{\preq}{\simeq_K}
\nc{\lb}{[\,,\,]}
\nc{\ds}{\displaystyle}

\nc{\ben}{\begin{enumerate}}\nc{\een}{\end{enumerate}}

\nc{\He}{\operatorname{Hess}} \nc{\ad}{\operatorname{ad}}
\nc{\Ad}{\operatorname{Ad}} \nc{\rank}{\operatorname{rank}}
\nc{\Irr}{\operatorname{Irr}} \nc{\End}{\operatorname{End}}
\nc{\Aut}{\operatorname{Aut}} \nc{\Inn}{\operatorname{Inn}}
\nc{\Der}{\operatorname{Der}} \nc{\Ker}{\operatorname{Ker}}
\nc{\Iso}{\operatorname{I}} \nc{\Diff}{\operatorname{Diff}}
\nc{\Lie}{\operatorname{L}} \nc{\tr}{\operatorname{tr}}

\nc{\dQ}{\operatorname{\mbox{deg}_\mathbb{Q}}}
\nc{\degr}{\operatorname{deg}} \nc{\sen}{\operatorname{sen}}
\nc{\modu}{\operatorname{mod}} \nc{\Ric}{\operatorname{Ric}}
\nc{\sym}{\operatorname{sym}} \nc{\sca}{\operatorname{sc}}
\nc{\scalar}{{\sf s}} \nc{\grad}{\operatorname{grad}}
\nc{\ricci}{\operatorname{ric}} \nc{\Rin}{\operatorname{M}}
\nc{\Le}{\operatorname{L}}
\nc{\level}{\operatorname{level}} \nc{\rad}{\operatorname{r}}
\nc{\abel}{\operatorname{ab}} \nc{\Pf}{\operatorname{Pf}}

\newtheorem{theorem}{Theorem}[section]
\newtheorem{proposition}[theorem]{Proposition}

\newtheorem{remark}[theorem]{Remark}

\newtheorem{example}[theorem]{Example}

\title{Examples of Anosov Lie algebras}

\author[Mainkar,\quad Will]{Meera G. Mainkar,\quad \dag Cynthia E. Will}

\address{School of Mathematics, Tata Institute of Fundamental Research, 
 Mumbai, India// FaMAF and CIEM, Universidad Nacional de C\'ordoba, Haya de la Torre s/n,
5000 C\'ordoba, Argentina } \email{meera@math.tifr.res.in,
cwill@mate.uncor.edu}

\thanks{
\dag Supported by CONICET, SeCyT and FONCyT (Argentina). }

\begin{document}

\begin{abstract}
We construct new families of  examples of (real) Anosov Lie
algebras starting with algebraic units. We also give examples of
indecomposable  Anosov Lie algebras (not a direct sum of proper
Lie ideals) of dimension $13$ and $16$, and we conclude that for
every $n \geq 6$ with $n \neq 7$ there exists an indecomposable
Anosov Lie algebra of dimension $n$.
\end{abstract}

\maketitle

\section{Introduction}\label{intro}
A diffeomorphism $f$ of a compact differentiable manifold $M$ is
called {\it Anosov} if it has a global hyperbolic behavior, i.e.
the tangent bundle $\tang M$ admits a continuous invariant
splitting $\tang M=E^+\oplus E^-$ such that $\dif f$ expands $E^+$
and contracts $E^-$ exponentially. This kind of diffeomorphism
plays an important and beautiful role in dynamics since they give
examples of dynamical systems with very nice properties, and it is
then a natural problem to understand which are the manifolds
supporting them (see \cite{SS}).

Up to now, the only known examples are hyperbolic automorphisms on
infranilmanifolds (manifolds finitely covered by  nilmanifold)
which are called {\it Anosov automorphisms}. Moreover, it is
conjectured that any Anosov diffeomorphism is topologically
conjugate to an Anosov automorphism of a infranilmanifold (see
\cite{Mrg}). The conjecture is known to be true in many particular
cases, for example, J. Franks \cite{F} and A. Manning \cite{Man}
proved it for Anosov diffeomorphisms on infranilmanifolds
themselve.

 We will say that an $n$-dimensional rational Lie algebra is {\it Anosov} if it
 admits a hyperbolic automorphism $\tau$ (i.e. none of the
 eigenvalues of $\tau$ are of modulus $1$) such that
 $[\tau]_{\beta}\in\Gl_n(\ZZ)$ for some basis $\beta$ of $\ngo$,
 where $[\tau]_{\beta}$ denotes the matrix of $\tau$ with respect to
 $\beta$. We say that a real Lie algebra is {\it Anosov} if it
 admits a rational form which is Anosov. It is easy to observe that
 a real Lie algebra $\ngo$ is Anosov if and only if it admits a hyperbolic
 automorphism $\tau$ such that
 $[\tau]_{\beta}\in\Gl_n(\ZZ)$
 for some $\ZZ$-basis $\beta$ of $\ngo$ (i.e.
with integer structure constants).

It is well known that any Anosov Lie algebra is necessarily
nilpotent, and it is easy to see that the classification of
nilmanifolds which admit an Anosov automorphism is essentially
equivalent to that of Anosov Lie algebras (see \cite{L,Dn,Ito,D}).

Therefore, if one is interested in finding those Lie groups which
are simply connected covers of Anosov infranilmanifold, then the
objects to find are real nilpotent Lie algebras $\ngo$ supporting
an Anosov automorphism.

Concerning the known examples, beside the case of free nilpotent
Lie algebras (see \cite{Dn}), there were only sporadic examples of
Anosov Lie algebras before \cite{L}, where  it is proved that
$\widetilde{\ngo} = \ngo+...+\ngo$ is a real Anosov Lie algebra
for any graded Lie algebra $\ngo$ admiting a rational form. Also,
in \cite{D-M} other kind of examples are given in the context of
certain two-step nilpotent Lie algebras attached to graphs. In
this way, there are in the literature examples of nonabelian
Anosov real Lie algebra for each dimension $n \geq 6$ but $7$ and
$13$. Moreover, in \cite{D-M} the existence of indecomposable
$n$-dimensional $2$-step Anosov Lie algebra
 is proved for  $n \geq 6,$ except for  $n = 7, 9, 12, 13, 16$. We recall
 that a Lie algebra is said to be  {\it indecomposable} if it can
 not be expressed as a direct sum of proper Lie ideals.
 It is known that there is no $7$-dimensional Anosov
Lie Algebra \cite{lw} and for $n = 9, 12$, there exists an
indecomposable Anosov Lie algebra of dimension $n$ (see \cite{L}).
In fact, \cite{L} gives a family of indecomposable  $(3r +
3)$-dimensional Anosov Lie algebras, $r \geq 1$.

In this paper we will give explicit families of examples of Anosov
(real) Lie algebras to illustrate a general procedure to construct
Anosov Lie algebras and, as an application we will give an
indecomposable $13$-dimensional Anosov Lie algebra. In fact, for
each pair of algebraic integers $\l,\;\mu$ of degree $p$ and $q$
respectively which satisfy \ben \item they are units, \item if we
denote by $\{\l=\l_1,\dots,\l_p\}$ and $\{\mu=\mu_1,\dots,\mu_q\}$
the conjugates to $\l$ and $\mu$ respectively, then $|\l_i|\ne
1\ne |\mu_j|$, and  \item $|\l_i\l_j| \ne 1,$ \een we will exhibit
a type $(p,q)$ Anosov Lie algebra. This first construction is
quite easy to extend and we are able to show examples of $3$-step
(and in fact of $k$-step) Anosov Lie algebras, and also in the
special case of $p=2$ we give another example of type $(3q,q+2)$
for any $q$.

 Finally, we also
give an example of an indecomposable $16$-dimensional Anosov Lie
algebra, which allows us to conclude that for $n \geq 6, n \neq 7$
there exists an indecomposable $n$-dimensional Anosov Lie algebra.

\vs

 {\bf Acknowledgements}: We would like to thank Dr. Miatello for useful comments and
 Dr. Lauret for
his invaluable help. The first author would also want to thank
TWAS for supporting her to stay in C\'ordoba for three month and
to CIEM for the hospitality during that stay.

\section{Examples}\label{ex}

Given a nilpotent Lie algebra $\ngo$, we call the {\it type} of
$\ngo$ to the $r$-tuple $(n_1,...,n_r)$, where
$n_i=\dim{C^{i-1}(\ngo)/C^i(\ngo)}$ and $C^i(\ngo)$ is the central
descending series. It is proven in \cite{lw} that if  $\ngo$ is a
real Anosov Lie algebra of type $(n_1,...,n_r)$, then there exist
a hyperbolic $A\in\Aut(\ngo)$ such that
\begin{itemize}
\item[(i)] $A\ngo_i=\ngo_i$ for all $i=1,...,r$,

\item[(ii)] $A$ is semisimple (in particular $A$ is diagonalizable
over $\CC$),

\item[(iii)] For each $i$, there exists a basis $\beta_i$ of
$\ngo_i$ such that $[A_i]_{\beta_i}\in\Sl_{n_i}(\ZZ),$ where
$n_i=\dim{\ngo_i}$ and $A_i=A|_{\ngo_i}$.

\end{itemize}

It is important to mention that the existence of an Anosov
automorphism is a really strong condition on an infranilmanifold
and also in a Lie algebra, and therefore, our approach is to start
with a hyperbolic automorphism.

In this context, to show an example of an Anosov Lie algebra, we
are going to construct a complex Lie algebra (to be able to work
with eigenvalues) in such a way that it admits a hyperbolic
automorphism $A$ such that $[A]_{\beta}\in\Gl_n(\ZZ)$ for some
$\ZZ$-basis $\beta$ of $\ngo.$

We begin by noting that if $\l$ and $\mu$ are algebraic units of
degree $p$ and $q$ respectively, and we denote by $\{\l=\l_1,
\dots, \l_p\}$ and $\{\mu=\mu_1, \dots, \mu_q\}$ the sets of
conjugates of $\l$ and $\mu$ over $\QQ$ respectively, it is not
hard to see that $\{\l_i\mu_j\}$ are also algebraic units and
moreover, the matrix
$\left[\begin{smallmatrix}\l_1\mu_1&&\\&\ddots&\\&&\l_p\mu_q\end{smallmatrix}\right]$
is conjugated to a matrix in $\Gl_{pq}(\mathbb{Z})$ with
determinant $\pm 1$.

Bearing this in mind, for each pair of non negative integers
$p\ne q$ , we take the Lie algebra $\ngo$ with basis $\beta=\{X_1,
\dots, X_{pq}, Y_1, \dots, Y_p,Z_1, \dots, Z_q\}$ and  Lie bracket
among  given by:
\begin{equation}\label{npq}
[X_{ip+j},Y_j]=Z_{i+1} \qquad 0\le i <q,\;1\le j \le p.
\end{equation}
It is clear that $\ngo$ is a two-step nilpotent Lie algebra, a basis
of $\ngo_1$ is  $\{X_i, Y_j : 1\le i\le pq, \;1\le j\le p\}$ and
$\{Z_k : \; 1 \le k \le q\}$ is a basis for $\ngo_2.$  Now, let $A$
be an automorphism such that $[A]_{\beta} =\left[\begin{smallmatrix}
A_1&
\\& A_2
\end{smallmatrix}\right]$ where
$$\begin{array}{lcl}
A_1=\left[\begin{smallmatrix}
  \l_1\mu_1 &        &          &          &       &          &          &       & \\
            & \ddots &          &          &       &          &          &       & \\
            &        & \l_p\mu_1&          &       &          &          &       & \\
            &        &          & \l_1\mu_2&       &          &          &       & \\
            &        &          &          & \ddots&          &          &       & \\
            &        &          &          &       & \l_p\mu_q&          &       & \\
            &        &          &          &       &          & \l_1^{-1}&       & \\
            &        &          &          &       &          &          &\ddots & \\
            &        &          &          &       &          &          &       & \l_p^{-1}
\end{smallmatrix}\right] & \mbox{and} &
A_2=\left[ \begin{array}{ccc} \mu_1 &&\\ &\ddots&\\&&\mu_q
\end{array}\right].
\end{array}$$

We note that $\beta$ is a basis of eigenvectors of $A$. Also, if
we take $\l$ and $\mu$ as above and such that $|\l_i| \neq 1,
|\mu_j| \neq 1$ and $|\l_i\mu_j| \ne 1$ for all $i,j$ then $A$ is
a hyperbolic automorphism.

In what follows, we are going to show that $\ngo$ is an Anosov Lie
algebra by constructing a $\mathbb{Z}$-basis of $\ngo$ preserved
by $A.$ In order to make the calculation more clear we will make a
small change in the notation. Let $X_{(i,j)}$ be the eigenvector
of $A$ corresponding to the eigenvalue $\l_i\mu_j$. Note that this
is only a reordering of the $\{X_i\}.$ In fact, $X_{(i,j)}=
X_{(j-1)p+i},$ and therefore we may say that
$\beta=\{X_{(i,j)},Y_k,Z_l: 1\le i,k \le p,\; 1\le j,l \le q \} $
and (\ref{npq}) is now given by
\begin{equation}\label{npq1}
[X_{(i,j)},Y_i]=Z_j.
\end{equation}
Let $\beta'=\{\mathcal{X}_{(k,l)},\mathcal{Y}_r ,\mathcal{Z}_s :
0\le i,k < p,\; 0\le j,l < q \} $ be the new basis of $\ngo$ given
by
$$\begin{array}{ll}
\mathcal{X}_{(k,l)}=\ds{\sum_{i=1}^{p}\sum_{j=1}^{q}}
\l_i^{k}\mu_j^{l}\, X_{(i,j)} & 0\le k < p,\; 0\le l < q,  \\
& \\
\mathcal{Y}_r=\ds{\sum_{k=1}^{p}} \l_k^{-r} Y_k &  0\le r < p, \\
&\\
 \mathcal{Z}_s=\ds{\sum_{l=1}^{q}} \mu_l^s Z_l & 0\le s < q.
\end{array}$$

To see that this is actually a basis of $\ngo$, it is enough to
check that the sets $\{\mathcal{X}_{(k,l)}\}$, $\{\mathcal{Y}_r\}$
and $\{\mathcal{Z}_s\}$ are linearly independent over $\QQ$. Since
all the calculations are similar, we are only going to show how to
proceed with $\{\mathcal{X}_{(k,l)}\}$. Suppose $a_{kl}'$s $\in
\QQ$ such that
$$\begin{array}{lcl}
0&=&\ds{\sum_{k=0}^{p-1}\sum_{l=0}^{q-1}}\;
a_{kl}\mathcal{X}_{(k,l)}\\&&\\
&=&\ds{\sum_{k=0}^{p-1}\sum_{l=0}^{q-1}} \,
a_{kl}\left(\ds{\sum_{i=1}^{p}\sum_{j=1}^{q}} \l_i^k\mu_j^l \,
X_{(i,j)}\right)\\ &&\\
&=&\ds{\sum_{i=1}^{p}\sum_{j=1}^{q}}\left(\ds{\sum_{k=0}^{p-1}\sum_{l=0}^{q-1}}
a_{kl}\l_i^k\mu_j^l \right) \, X_{(i,j)}.
\end{array}
$$
Hence, for $1\le i \le p,\; 1\le j \le q$ we have that
$$
0=\ds{\sum_{k=0}^{p-1}}\left(\ds{\sum_{l=0}^{q-1}}a_{kl}\mu_j^l
\right)\,\l_i^k.
$$
This can be seen, for each $1\le j \le q$ fixed, as a polynomial
in $\l_i.$ This polynomial has degree $p-1$ and it vanish on each
one of the $\l_i,$ so by our choice of $\l$, it has $p$ different
roots and therefore is identically zero. Hence for $1\le j \le q$
we have that its coefficients are zero. That is, for each $0\le k
< p$
$$0=\ds{\sum_{l=0}^{q-1}}a_{kl}\mu_j^l $$
which is again a polynomial in $\mu_j$ of degree $q-1$ with $q$
different roots,  and therefore we can conclude that $a_{kl}=0$
for all $k,l$ as we wanted to show.

If $x^p+a_{p-1}x^{p-1}+\dots+ a_0$ and $x^p+b_{q-1}x^{q-1}+\dots+
b_0$ are the minimal polynomial of $\l^{-1}$ and $\mu$
respectively, it is not hard to check that
$$\begin{array}{lcl}
A \mathcal{Y}_r = \left\{
 \begin{array}{ll} \mathcal{Y}_{r+1}
 \; &  r < p-1,\\ &  \\
 -\ds{\sum_{j=0}^{p-1}} a_j \mathcal{Y}_j \; &  r = p-1,
 \end{array}\right. & \quad &
 A \mathcal{Z}_s = \left\{
 \begin{array}{ll} \mathcal{Z}_{s+1} &  s < q-1,\\ &  \\
 -\ds{\sum_{l=0}^{q-1}} b_l \mathcal{Z}_l &  s=
 q-1,\end{array}\right.
 \end{array}
 $$
  Note that $a_i$ and $b_j$ are all integer numbers.

Concerning  $\mathcal{X}_{(k,l)},$ by the definition we have that
for each $i,j$,
$$A (\l_i^k\mu_j^l\,X_{(i,j)})= \l_i^{k+1}\mu_j^{l+1}\,X_{(i,j)},$$
and therefore, for $k < p-1$ and $l < q-1$, $A\mathcal{X}_{(k,l)}=
\mathcal{X}_{(k+1,l+1)}.$ In the same line of the calculation done
above, we have that
$$ \begin{array}{l}
A\mathcal{X}_{(k,l)}= \left\{
 \begin{array}{ll} -\ds{\sum_{k=1}^{p-1}}c_k\,
\mathcal{X}_{(k,l+1)}&  k=p-1,\; l < q-1,\\
 -\ds{\sum_{l=1}^{q-1}}b_l\,
\mathcal{X}_{(k+1,l)}&  k < p-1,\; l=q-1\\
  \end{array}\right. \\
A\mathcal{X}_{(p-1,q-1)}= -\ds{\sum_{k=1}^{p-1}}
\sum_{l=1}^{q-1}c_k b_l\,\mathcal{X}_{(k,l)}
\end{array}$$

where $c_j \in \ZZ$ are the coefficients of the minimal polynomial
of $\l.$

On the other hand, to see that the Lie bracket of any two elements
of $\beta'$ is a linear combination of elements of $\beta'$ with
integer coefficients, it is enough to check it for
$[\mathcal{X}_{(k,l)},\mathcal{Y}_r].$ Using (\ref{npq1}) we have
that
\begin{equation}\label{corchete}
\begin{array}{lcl}
[\mathcal{X}_{(k,l)},\mathcal{Y}_r]&=&
\ds{\sum_{i=1}^{p}\sum_{j=1}^{q}}\l_i^{k-r}
\mu_j^l\,[X_{(i,j)},Y_i]\\&&\\ &=& \left(\ds{\sum_{i=1}^{p}}
\l_i^{k-r}\right)\left( \ds{\sum_{j=1}^{q}} \mu_j^l\,Z_j\right)\\& &\\
&=& M(k,l) \mathcal{Z}_l.
\end{array}
\end{equation}
Here $M(k,l)= \tr A_\l^{k-l}$ where
$A_\l=\left[\begin{smallmatrix}\l_1&&\\&\ddots&\\&&\l_p\end{smallmatrix}\right].$
Due to our choice of $\l,$ $A_\l$ is conjugated to a matrix in
$\Gl_p(\mathbb{Z})$ and therefore so is $A_\l^m$ for any $m \in
\mathbb{N}.$  Hence $M(k,l)$ is an integer number for any $k,l$ as
we wanted to show.

\begin{remark}\label{rem1}{\rm Note that the Lie algebra $\ngo$ we have
constructed does not depend on the algebraic numbers $\l$ and
$\mu$, it only depends on $p$ and $q$, and moreover it is easy to
see (by looking at the dimension of the center for example) that
the Lie algebra associated to $(p,q)$ is not isomorphic to the one
corresponding to $(q,p)$ unless $p=q$. We have obtained in this
way two non isomorphic Anosov Lie algebra of dimension $n$ for all
$n=p.q+p+q$ for any non negative integers $p,q$. It is easy to
check that for $p=2=q,$ we obtain the two step nilpotent Lie
algebra $\ggo,$ of type $(6,2)$ given in \cite{lw}.

Concerning the existence of algebraic numbers as we need, we refer
to \cite{Li}.

Finally we would like to point out, for further use, that $\ngo$
can be viewed as $V_0 \oplus V_1 \oplus Z$ where $V_0$ is the
subspace generated by the $\{X_{(i,j)}\}$, $V_1$ is the one
spanned by the $\{Y_k\}$ and $Z$ is the center. In this setup,
$V_0$ acts on $V_1\oplus Z$, as it is stated in (\ref{npq1}).}
\end{remark}

\begin{example}{\rm
As a new example, we can carry out the calculations for $p=2,\;
q=3$ to obtain the $11$-dimensional Lie algebra with basis
$$\beta=\{ X_1,\dots,X_6,Y_1,Y_2,Y_3,Z_1,Z_2\}$$ and Lie bracket
among them given by
\begin{equation}\label{n23}
\begin{array}{lcl}
[X_1,Y_1]=Z_1 & [X_3,Y_2]=Z_1 & [X_5,Y_3]=Z_1 \\
& &\\

[X_2,Y_1]=Z_2 & [X_4,Y_2]=Z_2 & [X_6,Y_3]=Z_2. \\
\end{array}
\end{equation}
The hyperbolic automorphism $A$ is given by
$[A]_\beta=\left[\begin{smallmatrix}A_0& &\\&A_1&\\&
&A_2\end{smallmatrix}\right]$ where
$$
\begin{array}{lcl}
A_0=\left[\begin{smallmatrix}\l\mu_1 & & & & & \\
                               & \l^{-1}\mu_1 & & & & \\
                                && \l\mu_2 &&& \\
                               &&& \l^{-1} \mu_2 && \\
                               &&&& \l\mu_3 & \\
                               &&&& & \l^{-1}\mu_3  \end{smallmatrix}\right],&
A_1=\left[\begin{smallmatrix}\mu_1^{-1}&&\\&\mu_2^{-1}&\\&&\mu_3^{-1}\end{smallmatrix}\right],
& A_2
\left[\begin{smallmatrix}\l&\\&\l^{-1}\end{smallmatrix}\right].
\end{array}$$

In this case we have obtained a Lie algebra of type $(9,2),$ and
note that for $p=3,\; q=2$ we obtain a Lie algebra of type
$(8,3).$ We would like to point out that here and in general, we
can add non zero constant to the Lie brackets in (\ref{n23}) but
it is easy to see that this leads to isomorphic Lie algebras.}
\end{example}

Once we have stated the general picture, let us consider an
analogous procedure by starting from two algebraic units $\l$ and
$\mu.$ In this case, by following essentially the same procedure as
above, we can construct a two step nilpotent Anosov Lie algebra of
type $(p+q,pq)$, where the eigenvalues of the corresponding $A_1$
are the conjugated numbers to $\l$ and $\mu$ and the ones
corresponding to $A_2$ are all the products among them,
$\{\l_i\mu_j\}.$ It is not hard to see that this algebra is
isomorphic to the one associated to a bipartite graph $(p,q)$ which
is proved to be Anosov in \cite{D-M}. As in this case, a lot of
changes can be made to this procedure to obtain a variety of new
examples. Among them we are now going to mention a few more, and
since the proofs are essentially the same, they will be omitted.

\vspace{.1cm}

{\bf Example 1.} One can start be taking three algebraic units as
above $\l,\mu$ and  $\nu$ of degree $p,q,$ and $r$ respectively,
such that the conjugate numbers to each of them satisfies  $|\l_i
\mu_j| \ne 1$, $|\l_i \nu_k| \ne 1$ and $|\nu_k \mu_j| \ne 1.$ It
is not hard to see that we can proceed analogously considering the
pair $\l\mu$ and $\nu,$ $(pq,r)$ in spite of which is the degree
of $\l\mu.$ In fact, in the proof of the linear independence of
the new basis, and also in (\ref{corchete}), we only use the fact
that we are adding over all the conjugated numbers to $\l$ and
$\mu.$ Following the lines of the above procedure, we then obtain
an Anosov Lie algebra,
 $\ngo_{(pq,r)}$ of type $(pqr+r,pq).$
Moreover, once we have stated this, it is clear that it is
 also true for $\l\nu$ and $\mu,$ $(pr,q)$ and in this
case, our procedure leads to a Lie algebra of type $(prq+q,pr)$.

Now, it is clear that in each of this algebras, $\ngo_{(pq,r)}$
and  $\ngo_{(pr,q)},$ the eigenvalues of the associated
automorphism corresponding to ${X}_{(k,l)}$ are the same, that is
${\l_i\mu_j\nu_s}$ for some ${i,j,s}$. Therefore, the
corresponding subspaces $V_0$ can be identifying (see remark
\ref{rem1}). In this case, it is easy to see that a new algebra
can be constructed from this two by identifying the $V_0$.
Explicitly, if $\ngo_{(pq,r)}=V_0\oplus V_1 \oplus Z_1$ and
$\ngo_{(pr,q)}=V_0\oplus V_2 \oplus Z_2$, let $\ngo$ be the Lie
algebra with vector space $\ngo=(V_0\oplus V_1\oplus V_2) \oplus
(Z_1 \oplus Z_2),$ where the action is as before: $[V_0,V_i]
\subset Z_i,$ $i=1,2$. This is a two step nilpotent Lie algebra of
type $(pqr+q+r,pr+pq)$. In this framework, there is a natural way
to define an automorphism in $\ngo,$ using the ones in
$\ngo_{(pq,r)}$ and $\ngo_{(pr,q)}:$
$[A]_\beta=\left[\begin{smallmatrix}A_1&
\\&A_2\end{smallmatrix}\right]$ where
$$\begin{array}{ll}
A_1=\left[\begin{smallmatrix}\l_1\mu_1\nu_1 & & & & &&&& \\
                               & \ddots && & & &&& \\
                               && \l_p\mu_q\nu_r &&&&&&\\
                                &&& \mu_1^{-1} &&&&& \\
                               &&&&\ddots &&&& \\
                               &&&&& \mu_q^{-1}&&& \\
                               &&&&& & \nu_1^{-1}&&\\
                               &&&&&&&\ddots&\\
                               &&&&&&&&\nu_r^{-1}  \end{smallmatrix}\right],&
A_2=\left[\begin{smallmatrix}\l_1\nu_1&&&&&\\&\ddots&&&&\\&&\l_p\nu_r&&&\\&&&\l_1\mu_1&&\\
&&&&\ddots&\\ &&&&&\l_p\mu_q \end{smallmatrix}\right].
\end{array}$$
It is easy to check that due to our choice of $\l$, $\mu$ and $\nu$, it is hyperbolic.
On the other hand, note
that in both cases the lattice we have constructed in $V_0,$ ${\mathcal{X}_{(k,l)},}$
are the same (is just a
matter of notation) and therefore as with the automorphism, the natural extension of the
lattice we have in
$\ngo_{(pq,r)}$ and $\ngo_{(pr,q)},$ is a $\mathbb{Z}$ basis preserved by $A$ and therefore
$\ngo$ is an Anosov
Lie Algebra.

In this way we obtain two step Anosov Lie algebras of dimension
$n=pqr+pq+pr+q+r$ for any $p,q,r$. Distinguishing them by the
type, it is easy to see that in general, if $p\ne q \ne r$ the Lie
algebra one obtains by interchanging the roll of $p,\,q$ and $r$
are not isomorphic. The smallest one we can construct corresponds
to $p=q=r=2$ is $18$ dimensional and its type is $(10,8)$.

It is not hard to see that this procedure extends in a natural way
to consider $k$ algebraic units to obtain a $2$-step Anosov Lie
algebra.

\vspace{.1cm}

{\bf Example 2.} Now we are going to show how to use the procedure
to construct three-step Anosov Lie algebras. As before we take
algebraic numbers as above $\l,\mu$ and  $\nu$ of degree $p,\;q,$
and $r$ respectively.

In this case we have in mind
$A=\left[\begin{smallmatrix}A_1&&\\&A_2&\\&&A_3\end{smallmatrix}\right]$,
where $A_1$ and $A_2$ are as in the previous example, and
$A_3=\left[\begin{smallmatrix}\l_1&&\\&\ddots&\\&&\l_p\end{smallmatrix}\right].$
As before, we are going to make a small change in the notation in
order to be consistent with the eigenvalues. Let $\ngo$ be the Lie
algebra with basis
$$\beta=\left\{ X_{(i,j,k)}, Y_j,Z_k,V_{(i,k)},W_{(i,j)},U_i : 1\le i\le
p,\;1\le j\le q,\;1\le k\le r \right\}$$ and the Lie bracket among
them given by
$$\begin{array}{lcl} [X_{(i,j,k)},Y_m]=
\delta_{(j,m)}V_{(i,k)}&\quad &
[X_{(i,j,k)},Z_n]= \delta_{(k,n)}W_{(i,j)} \\
& & \\

[Z_n,V_{(i,k)}]=\delta_{(n,k)}U_i &\quad &

[Y_m,W_{(i,j)}]=\delta_{(m,j)}U_i.

 \end{array}$$

It is easy to see that $\ngo$ is a three-step nilpotent Lie
algebra, that is, it satisfies Jacobi identities, and the type of
$\ngo$ is $(pqr+q+r,pr+pq,p).$ Let $A$ denote the linear
transformation of $\ngo$ such that
$[A]_\beta=\left[\begin{smallmatrix}A_1&&\\&A_2&\\&&A_3\end{smallmatrix}\right]$
hence, $A$ is a hyperbolic automorphism of $\ngo$ and $\beta$ is a
basis of eigenvectors of $A$. To construct a $\mathbb{Z}$-basis,
we proceed similarly as before:
$$\begin{array}{ll}
\mathcal{X}_{(m,l,s)}=\ds{\sum_{i=1}^{p}\sum_{j=1}^{q}\sum_{k=1}^r}
\l_i^{m}\mu_j^{l}\nu_k^s\, X_{(i,j,k)} & 0\le m < p,\; 0\le l < q,\;0\le k < r,   \\
& \\
\mathcal{Y}_l=\ds{\sum_{k=1}^{q}} \mu_k^{-l}\; Y_k &  0\le l < q,
\\ &\\
 \mathcal{Z}_s=\ds{\sum_{l=1}^{r}} \nu_l^{-s}\; Z_l & 0\le s < r,\\
 &\\
 \mathcal{V}_{(m,s)}=\ds{\sum_{i=1}^{p}\sum_{k=1}^{r}}
\l_i^{m}\nu_k^{s}\, V_{(i,k)} & 0\le m < p,\; 0\le s < r,\\
&\\ \mathcal{W}_{(m,l)}=\ds{\sum_{i=1}^{p}\sum_{j=1}^{q}}
\l_i^{m}\mu_j^{l}\, W_{(i,j)} & 0\le m < p,\; 0\le l < q,\\
&\\
 \mathcal{U}_i=\ds{\sum_{n=1}^{p}} \l_n^{i}\; U_n &  0\le i <
p.
\end{array}$$

Straightforward calculation, using the same techniques as before,
shows that this is a $\mathbb{Z}$-basis preserved by $A$, and
therefore $\ngo$ is an Anosov Lie algebra. The smallest example we
obtain in this way is a tree step nilpotent Lie algebra of
dimension  $20$ of type $(10,8,2)$.

It is not hard to see that if $\ngo=\ngo_1\oplus \ngo_r$  is a real
Anosov Lie algebra of type $(n_1,n_2,\dots,n_r)$  then $\ngo/\ngo_r$
is also an Anosov Lie algebra (see \cite{mw}). Note that in this
case, this fact is what we have showed in the previous construction.

Also, it is not hard to prove by induction, that this procedure
extends in a natural way to the case of consider $k$-algebraic
units, to obtain a $k$-step Anosov Lie algebra. Moreover, by the
above observation, we have all the quotients one has in between.

\vspace{.3cm}

{\bf Example 3.} As a last application of our procedure we are
going to consider the special case of $p=2.$ In this case, in
addition to $\ngo_{(2,q)},$ and $\ngo_{(q,2)},$ we can define
others Lie algebras, for example by adding Lie brackets among the
$\{X_{(i,j)}\}$. That is, we take
$$
\beta=\{ X_i, Y_k, Z_l:\; 1\le i\le 2q\;1\le k\le q,\; 1\le l\le
q+2 \}
$$
as a basis of $\ngo$ with the Lie bracket given by
\begin{equation}
\begin{array}{lcl}
[X_{iq+j},Y_j]= Z_{q+i+1} & \qquad & i=0,1:\, j=1\dots,q \\
&&\\

[X_j,X_{j+q}]= Z_j & & j=1\dots,q
 \end{array}
\end{equation}
This is a two-step nilpotent Lie algebra of type $(3q,q+2)$. Let
$A$ be the automorphism of $\ngo$ such that
$[A]_\beta=\left[\begin{smallmatrix}A_1&\\&A_2\end{smallmatrix}\right],$
where
$$\begin{array}{ll}
A_1=\left[\begin{smallmatrix}\l\mu_1 & & & & & \\
                               & \ddots && & &  \\
                               && \l^{-1}\mu_q &&&\\
                                &&& \mu_1^{-1} && \\
                               &&&&\ddots & \\
                               &&&&& \mu_q^{-1}
                                 \end{smallmatrix}\right],&
A_2=\left[\begin{smallmatrix}\mu_1^2&&&&\\&\ddots&&&\\&&\mu_q^2&&\\&&&\l&\\
&&&&\l^{-1} \end{smallmatrix}\right].
\end{array}$$
Concerning the lattice, we will take $\mathcal{Y}_k$ as before,
and let

\begin{equation}\label{ex3}
 \begin{array}{l} \mathcal{X}_i=\left\{
\begin{array}{ll}\ds{\sum_{k=1}^{q}}
\mu_k^{i}\,\left( X_k + X_{q+k}\right) &\qquad 0 \le i< q\\
\ds{\sum_{k=1}^{q}}
\mu_k^{i-q}\, \left( \l\, X_k + \l^{-1}\,X_{q+k}\right) &\qquad q \le i< 2q\\
\end{array}\right. \\ \\
\mathcal{Z}_l= \left\{
\begin{array}{ll} (\l^{-1}-\l)\ds{\sum_{k=1}^{q}}
\mu_k^{l}\, Z_k & 0 \le l< q\\
Z_q+Z_{q+1} & l=q\\ & \\
\l Z_q+ \l^{-1}Z_{q+1} & l=q+1.
\end{array}\right.
\end{array}
\end{equation}

 One can see that this is also a basis of $\ngo$ preserved by
$A$, and moreover one can check that
\begin{equation}
\begin{array}{lcl}
[\mathcal{X}_i,\mathcal{X}_j]= 0 &\quad& \;  i,j < q\;\;
\mbox{or}\; i,j \ge q
\\
&&\\

[\mathcal{X}_i,\mathcal{X}_{j}]= \mathcal{Z}_{i + j - q},&& 0\le i <q
 \mbox{ and } q \leq j < 2q.\\\
 \end{array}
\end{equation}
 Also, as in (\ref{corchete}), one can see that
 $$[\mathcal{X}_{iq+j},\mathcal{Y}_{k}]=
 N(j,k) \mathcal{Z}_{q+i},\qquad i = 0, 1 : 0 \leq j < q,$$
where $N(j,k)=\tr (A_{\mu}^{j-k})$, such that $A_\mu =
\left[\begin{smallmatrix}
 \mu_1&&&\\&\mu_2&&\\&&\ddots&\\&&&\mu_q \end{smallmatrix}\right]$,
  is an integer number and therefore, we
can conclude that $\ngo$ is an Anosov Lie algebra. Note that the
dimension of $\ngo$ is $n=4q+2$ and the type is $(3q,q+2).$ The
small one we can obtain corresponds to $q=2$, is of type $(6,4)$,
dimension $10$.
\begin{remark}\label{ex3rem}
{\rm Note that the subalgebra generated by the set $\{ X_i, Z_l:\;
1\le i\le 2q,\;1\le l\le q \}$ is an Anosov Lie subalgebra of
$\ngo$}
\end{remark}

\vspace{.3cm}

\section{$13$-dimensional example}

\vspace{.3cm}

This last example is rather different. We will take $p=2$, $q=3$
and we will split the basis in the center. So let $\l$ and $\mu$
two algebraic numbers of degree $2$ and $3$ respectively such that
$|\l_i\mu_j|\ne 1$ for all its conjugated numbers of $\l$ and
$\mu$. Then we take $\ngo$ the complex vector space with basis
$$
\beta=\{ X_1, \dots, X_6, Y_1,Y_2,Y_3, Z_1,Z_2, W_1,W_2 \}
$$
and define the Lie bracket among them by
\begin{equation}\label{13dim}
\begin{array}{lcl}
[X_1,Y_1]= Z_1 & \qquad& [X_4,Y_1]= W_1 \\
&&\\

[X_2,Y_2]= Z_2 & \qquad& [X_5,Y_2]= W_2 \\
&&\\

[X_3,Y_3]=- (Z_1+Z_2)& \qquad & [X_6,Y_3]= - (W_1+W_2).
 \end{array}
\end{equation}
This is a two-step nilpotent Lie algebra of type $(9,4)$. Let $A$
be the automorphism of $\ngo$ such that
$[A]_\beta=\left[\begin{smallmatrix}A_1&\\&A_2\end{smallmatrix}\right],$
where
$$\begin{array}{ll}
A_1=\left[\begin{smallmatrix}\l\mu_1 & & & & & \\
                               & \ddots && & &  \\
                               && \l^{-1}\mu_3 &&&\\
                                &&& \mu_1^{-1} && \\
                               &&&&\mu_2^{-1} & \\
                               &&&&& \mu_3^{-1}
                                 \end{smallmatrix}\right],&
A_2=\left[\begin{smallmatrix}\l&&&&\\&\l&&&\\&&\l^{-1}&&\\&&&&\l^{-1}
\end{smallmatrix}\right].
\end{array}$$
Concerning the lattice, we will take $\{\mathcal{X}_i\},\, 0\le i
\le 5$ and $\{\mathcal{Y}_k\},\, 0 \le k\le 2$ as in the previous
example, and let
$$ \begin{array}{ll}
\mathcal{Z}_l= \ds\sum_{i=1}^{2}\left( \mu_i^{l}-\mu_3^{l}
\right)\left(Z_i+W_i\right) & l=-1,1,  \\
&\\
 \mathcal{W}_l= \ds{\sum_{i=1}^{2}}\left(
\mu_i^{l}-\mu_3^{l} \right)\left(\l\,Z_i+\l^{-1}\,W_i\right) & l=-1,1.
\end{array}
$$

To see that this is a basis of $\ngo$, as we have pointed out
before, it is enough to check that each one of
$\{\mathcal{X}_i\},\,\{\mathcal{Y}_j\}$ and
$\{\mathcal{Z}_k\,\,\mathcal{W}_l\}$ are linearly independent
sets. We also note that the calculations for the first two sets
have been already done, and then we are only going to show how to
proceed in the center. Suppose then that
$$\begin{array}{rl}
0 = &a_{-1}\mathcal{Z}_{-1} + a_{1}\mathcal{Z}_{1} +
      b_{-1}\mathcal{W}_{-1} + b_{1}\mathcal{W}_{1}
          \\
  = & \ds{\sum_{i=1}^{2}}
\Big[\left( \mu_i^{-1}-\mu_3^{-1} \right)\left( a_{-1}+\l
b_{-1}\right) + \left( \mu_i^{1}-\mu_3^{1} \right)\left(
a_{1}+\l b_{1} \right)\Big] \, Z_i  \\
& + \ds{\sum_{i=1}^{2}}\Big[ \left( \mu_i^{-1}-\mu_3^{-1}
\right)\left( a_{-1}+\l^{-1} b_{-1}\right) + \left(
\mu_i^{1}-\mu_3^{1} \right)\left( a_{1}+\l^{-1} b_{1} \right)\Big]
\, W_i
\end{array}
$$
Hence, for $i =1,2$ we have that
\begin{equation}\label{poly0}
\begin{array}{l} 0= \left( \mu_i^{-1}-\mu_3^{-1}
\right)\left( a_{-1}+\l b_{-1}\right) + \left( \mu_i^{1}-\mu_3^{1}
\right)\left( a_{1}+\l b_{1} \right),
\\ \\
0=\left( \mu_i^{-1}-\mu_3^{-1} \right)\left( a_{-1}+\l^{-1}
b_{-1}\right) + \left( \mu_i^{1}-\mu_3^{1} \right)\left(
a_{1}+\l^{-1} b_{1} \right).
\end{array}
\end{equation}
If we denote by
$$P_{\l}(x)= x^{-1}\left( a_{-1}+\l
b_{-1}\right) + x\left( a_{1}+\l b_{1} \right),$$ then by the
first equation we have that $P_{\l}(\mu_i) =P_{\l}(\mu_j)= C$ for
$i,j=1,2,3$. Hence,
\begin{equation} \label{poly1} \left( a_{-1}+\l
b_{-1}\right) + x^2\left( a_{1}+\l b_{1} \right)-Cx,
\end{equation}

is a degree two polynomial annulated by each one of the $\mu_i.$ Since these are
three different algebraic
numbers, we have that (\ref{poly1}) is identically zero. In particular
\begin{equation} \label{poly2}
 0=\left( a_{-1}+\l b_{-1}\right) \quad \mbox{and} \quad 0=\left(
a_{1}+\l b_{1} \right).\end{equation}
We can do the same for the
second equation in \ref{poly0} and we will obtain
\begin{equation} \label{poly3}
 0=\left( a_{-1}+\l^{-1} b_{-1}\right) \quad \mbox{and} \quad 0=\left(
a_{1}+\l^{-1} b_{1} \right).\end{equation}

Finally, from (\ref{poly2}) and (\ref{poly3}) we can conclude that
$$ a_{-1}= b_{-1}= a_{1}=b_{1}=0$$
as was to be shown.

One can also see for $i=1$ or $2$, let say for simplicity $i=1,$
that
$$\begin{array}{lcl}
\mu_1^{2}-\mu_3^{2}&=&\left(\mu_1-\mu_3\right)\left(\mu_1+\mu_3\right)\\
&&\\
&=& \left(\mu_1-\mu_3\right)\left(n_1-(\mu_1\mu_3)^{-1}\right)\\
&&\\
&=& n_1\left(\mu_1-\mu_3\right)+\left(\mu_1^{-1}-\mu_3^{-1}\right)
\end{array}$$
where $n_j=\tr A_\mu^j$ is an integer number for all $j \in
\mathbb{Z}$ and we have also used that $\mu_2=(\mu_1\mu_3)^{-1}$.
Therefore,
$$ \ds\sum_{i=1}^{2}\left( \mu_i^{2}-\mu_3^{2}
\right)\left(Z_i+W_i\right)=n_1\mathcal{Z}_1+ \mathcal{Z}_{-1}.$$

In the same way, we also have that
$$ \ds\sum_{i=1}^{2}\left( \mu_i^{2}-\mu_3^{2}
\right)\left(\l Z_i+\l^{-1}W_i\right)=n_1\mathcal{W}_1+
\mathcal{W}_{-1}$$
 It is also easy to see that this is also valid
for $\mu_i^{-2}-\mu_3^{-2}$, that is we have similar formulas to
these ones for
$$ \ds\sum_{i=1}^{2}\left( \mu_i^{-2}-\mu_3^{-2}
\right)\left(\l^j Z_i+\l^{-j}W_i\right),$$ for $j=0,1$.

Also, as in the previous examples, it is not hard to see that this
is a basis of $\ngo$ preserved by $A$. With all this, one can
check that
$$
\begin{array}{ll}
[\mathcal{X}_0,\mathcal{Y}_0]= 0 &
\;[\mathcal{X}_3,\mathcal{Y}_0]= 0
\\
&\\

[\mathcal{X}_1,\mathcal{Y}_0]= \mathcal{Z}_1 &
[\mathcal{X}_4,\mathcal{Y}_0]= \mathcal{W}_1,\\
&\\

[\mathcal{X}_2,\mathcal{Y}_0]= n_1\mathcal{Z}_1+ \mathcal{Z}_{-1}&
[\mathcal{X}_5,\mathcal{Y}_0]= n_1\mathcal{W}_1+
\mathcal{W}_{-1}\\
&\\

[\mathcal{X}_0,\mathcal{Y}_1]= \mathcal{Z}_{-1}&
[\mathcal{X}_0,\mathcal{Y}_2]= n_{-1}\mathcal{Z}_1+
\mathcal{Z}_{1}.
 \end{array}
$$

Using this and the fact that $A$ is an automorphism, it is easy to prove ,  that this is a $\mathbb{Z}$-basis of
$\ngo$. For example,
$$
\begin{array}{ll}
[\mathcal{X}_2,\mathcal{Y}_3]&=[A\mathcal{X}_1,A\mathcal{Y}_2]\\&\\
& = A[A\mathcal{X}_0,A\mathcal{Y}_1]\\&\\
 &=A(A\mathcal{Z}_{-1})\\
 &\\
&=A(\mathcal{W}_{-1})=-a\mathcal{W}_{-1}-\mathcal{Z}_{-1}
\end{array}$$
where $x^2+ax+1$ is the minimal polynomial of $\l$. Therefore, we
can conclude that this is a Anosov Lie algebra, as desired. In the
following we will prove by using similar arguments as in Lemma
$6.6$ of \cite{D-M}, that it is also indecomposable.

\begin{proposition}
 $\ngo$, defined by the relations given by (\ref{13dim}), is indecomposable.
\end{proposition}

\begin{proof} Suppose on the contrary that $\ngo = \ngo_1 \oplus \ngo_2$
is the sum of two nontrivial ideals of $\ngo.$ By definition, we
have that $\ngo = V \oplus W$
 where $V$ is the subspace
 spanned by the set  $S = \{X_1, \ldots, X_6,
Y_1, Y_2, Y_3\},$ and $W$ is the subspace  of $\ngo$ spanned by
$\{Z_1, Z_2, W_1, W_3 \}$. Let $p : \ngo \rightarrow V$ be the
projection onto $V$ with respect to this decomposition $\ngo = V
\oplus W$ and let $ V_1 = p(\ngo_1) $ and  $V_2 = p(\ngo_2)$. For
$i = 1, 2$ let
$$S_i = \{v \in S : \; a_v v + \ds{\sum_{v' \in S, v' \neq v}}
 a_{v'} v' \in V_i \; \mbox{for some } a_v \neq 0\}.$$

Then $S = S_1 \cup S_2,$ and since $\ngo_1$ and $\ngo_2$ are
nontrivial ideals, it is easy to see that $S_1$ and $S_2$ are
nonempty sets. Moreover, we have that $S_2 \setminus S_1$ is
empty. In fact, if  $S_2 \setminus S_1$ is nonempty, we can either
have that for all $v \in S_2 \setminus S_1$ and $v' \in S_1$, $[v,
v']$ is zero, or there exists $v \in S_2 \setminus S_1$ and $ v'
\in S_1$ such that $[v, v']$ is nonzero.

In the first situation, as $S = S_1 \cup S_2$, we may assume that
$Y_1$ is contained in $S_1$.
 Then there exists nonzero $a$ such that $a Y_1
+ x \in V_1$ where $x$ is contained in the span of $S \setminus
\{Y_1\}$. Since $\ngo_1$ is an ideal $[a Y_1 + x, X_1]$ and $[a
Y_1 + x, X_4]$ are contained in $\ngo_1$. This means $Z_1, W_1 \in
\ngo_1$.  We notice that not all $Y_i'$s are contained in $S_1$
because if all $Y_i'$s are contained in $S_1$ then (by our
assumption) all $X_j'$s are contained in $S_1$ and then $S_2
\setminus S_1$ is empty.  Now either $Y_2 \in S_1$ or $Y_2 \in S_2
\setminus S_1$. If $Y_2 \in S_1$ (similar argument works for
 the other case) then $Y_3$ must be contained in
$S_2 \setminus S_1$. In that case  $Z_2 $ and $W_2$ are contained
in $\ngo_1$ (by considering Lie brackets with $X_2$ and $X_5$),
  and similarly $Z_1 + Z_2$ and  $W_1 + W_2$ are contained in $\ngo_2$.
 This is a contradiction.

 On the other hand, if there exists $v \in  S_2
\setminus S_1$ and $ v' \in S_1$ such that $[v, v']$ is nonzero,
it is easy to see that if $v \in \{ Y_1, Y_2, Y_3 \}$ then $v' \in
\{X_1, \ldots X_6 \}$ and if $v \in \{X_1, \ldots X_6 \}$ then $v'
\in \{ Y_1, Y_2, Y_3 \}$. So, since it is entirely equivalent, we
can assume then that $v \in \{ Y_1, Y_2, Y_3 \}$ and $v = Y_1.$
Moreover, either $v' = X_1$ or $v' = X_4,$ and therefore, we may
assume that $v' = X_1$. From our definition of $S_i$, there exist
nonzero scalars $s$ and $t$ such that
 $s X_1 + x$ is contained in $V_1$ and $t Y_1 + y$ is contained in
$V_2$, where $x$ is in the subspace of $V$ spanned by $Y_1, Y_2,
Y_3, X_2, \ldots,  X_6$ and $y$ is in the subspace of $V$ spanned
by $Y_2, Y_3, X_1, X_2, \ldots, X_6$. Hence $[s X_1 + x, Y_1] \in
\ngo_1$
 and
 $[t Y_1 + y, X_1],  [t Y_1 + y, X_4] \in \ngo_2$
since $\ngo_1$ and $\ngo_2$ are  Lie ideals of $\ngo$. This
implies that
 $s Z_1 + s' W_1 \in \ngo_1$ where $s'$ is a scalar, and $Z_1, W_1 \in \ngo_2$.
 This is a contradiction since $s$ is nonzero and then we can conclude that
 $S_2 \setminus S_1$ is empty.

 Therefore, we have that $S = S_1$ and moreover we can see that $Z_1, Z_2, W_1, W_2$ are
contained in $\ngo_1$.  Hence $[\ngo_1 , \ngo_1] = [\ngo , \ngo]$
and from this, one has that $\ngo_2$ is in the center of $\ngo$.

On the other hand, it is easy to see that the center is equal to
$[\ngo , \ngo] = W$ and hence $\ngo_2$ is contained in $\ngo_1$,
contradicting our assumption that $\ngo_2$ is nontrivial. Hence
$\ngo$ can not be seen as a sum of two proper ideals as we wanted to
show. \end{proof}

\vspace{.3cm}

\section{$16$-dimensional example}

 Let $(S , E)$ denote the complete bipartite graph
 on a set  $S$ of $5$ elements
 partitioned into subsets $S_1$ and $S_2$ of $2$ and $3$ elements
 respectively. Following for example \cite{D-M} we can define from
 this (and any graph) a $2$-step nilpotent Lie
 algebra.
 Let $\nca = V \oplus W$
 denote the $2$-step nilpotent Lie algebra associated with $(S,
 E).$ We recall that in this case we obtain an Anosov Lie algebra of type $(5,6)$
 (see \cite{D-M}).
 Using this algebra we are going to construct a $16$-dimensional Anosov Lie
 algebra as follows.

 Let $\ngo$ be that Lie algebra with linear space $\ngo = V \oplus V \oplus W$ and Lie
 bracket defined by
 $[(x_1, x_2 , w) , (y_1 , y_2, w')] = [x_1 , y_1] + [x_2 , y_2]$ where
 $x_i$'s and $y_i$'s are vectors in $V$, $w, w' \in W$ and
  $[x_i , y_i]$ denotes the Lie bracket in $\nca$.
 To see that it is an Anosov Lie algebra, let us consider $\Phi$ the additive subgroup
  of $\ngo$ generated by
 the elements of the type $( v, 0, 0), (0, v', 0), (0, 0,
 [\gamma , \delta])$ where $v, v', \gamma, \delta \in S$.
 It is easy to see that
$\Phi$ is a $\ZZ$-subalgebra of $\ngo$ ( i.e. $\Phi$ is the set of
all $\ZZ$ linear combinations of the basis of $\ngo$ with integer
structure coefficients) and moreover, $\ngo$ admits a hyperbolic
automorphism $\tau$ such that $\tau (\Phi) = \Phi$. In fact, if
$\Phi'$  is a subgroup
 of $\nca$ generated by $ S \cup \{ [v , v'] : v , v' \in
 S \},$ $\nca$ admits a hyperbolic
 automorphism $\tau'$ such that $\tau'(\Phi') = \Phi'.$ (see \cite{D-M}, Theorem $1.1$).
 We take $\tau$ to be the natural extension of $\tau'$ to $\ngo$.
 Hence $\ngo$ is an Anosov Lie algebra.

\begin{proposition}
 $\ngo$, defined as above, is indecomposable.
\end{proposition}

\begin{proof}
 Let  $S_1 = \{ \alpha, \beta \}$ and $S_2 = \{ \gamma, \delta,
 \eta \}$.
 Suppose that $\ngo_1$ and $\ngo_2$ are two proper ideals of $\ngo$
 such that $\ngo = \ngo_1 \oplus \ngo_2$. As $[\ngo , \ngo] = [\ngo_1 , \ngo_1]
 \oplus [\ngo_2 , \ngo_2]$ and $[\ngo , \ngo]$ is $6$-dimensional, we may
  assume that dim $[\ngo_1 , \ngo_1] \leq 3$.
  Let $X = (v, v', w) \in \ngo_1$ where $v, v' \in V$ and $w \in W$.
  Let $v = \sum_{\zeta \in S}\,
  a_\zeta \zeta$. As $\ngo_1$ is an ideal, $[X , (\xi , 0 , 0)] \in \ngo_1$
 for all $\xi \in S$. Hence $a_\gamma \gamma \wedge \xi +
 a_\delta \delta \wedge \xi + a_\eta \eta \wedge \xi$
 and $a_\alpha \alpha \wedge \zeta + a_\beta \beta \wedge \zeta$ are
 contained in $[\ngo_1 , \ngo_1]$ for all $\xi \in S_1$ and $\zeta \in S_2$.
 As dim $[\ngo_1 , \ngo_1] \leq 3$, we see that either $a_\xi = 0$ for all
 $\xi \in S_1$ or $a_\zeta = 0$ for all $\zeta \in S_2$.
 Let $V_1$ (respectively $V_2$) denote the subspace  of $V$ spanned by
 $S_1$ (respectively $S_2$). Then by the above observation, $v \in V_1$
 or $v \in V_2$.

  Suppose
 $v \in V_1$. We will prove that
 $\ngo_1$ is contained in $ V_1 \oplus V \oplus W$.
 Suppose that  $v \neq 0$. Then the vectors
 $a_\alpha \alpha \wedge \zeta + a_\beta \beta \wedge \zeta \in
 [\ngo_1 , \ngo_1]$ for all $\zeta \in S_2$ are linearly independent.
  Hence  dim $[\ngo_1 , \ngo_1] = 3$.
 Suppose $(v_1, v_1', w_1) \in \ngo_1$ be such that  $v_1
 \in V_2$. Let $v_1 = a'_\gamma \gamma + a'_\delta \delta + a'_\eta
  \eta$. Now as dim $[\ngo_1 , \ngo_1] = 3$ and
 $a'_\gamma \gamma \wedge \alpha  +
 a'_\delta \delta \wedge \alpha + a'_\eta \eta \wedge \alpha $ and
 $a'_\gamma \gamma \wedge \beta  +
 a'_\delta \delta \wedge \beta + a'_\eta \eta \wedge \beta$ are
 contained in $[\ngo_1 , \ngo_1]$, $a'_\gamma =  a'_\delta =  a'_\eta =
 0$. Hence $v_1 = 0$. Thus we have proved that if $v \in V_1$, then
 $\ngo_1$ is contained in $ V_1 \oplus V \oplus W$. Similarly we prove
 that if $v \in V_2$, then $\ngo_1$ is contained in $V_2 \oplus V \oplus W$.
  Suppose $v \in V_2$ and $v \neq 0$.
  Let  $(v_1, v_1', w_1) \in \ngo_1$ be such that  $v_1
  \in V_1$, and write  $v_1 = a'_\alpha \alpha + a'_\beta \beta$.
 We note that  $a'_\alpha \alpha \wedge
 \zeta + a'_\beta \beta \wedge \zeta \in [\ngo_1 , \ngo_1]$ for all  $\zeta
 \in S_2$. If  the vectors
 $a'_\alpha \alpha \wedge
 \zeta + a'_\beta \beta \wedge \zeta$ are linearly independent for all
 $\zeta
 \in S_2$, then  dim $[\ngo_1 , \ngo_1] = 3$. This is  a
 contradiction as  $a_\gamma \gamma \wedge \alpha  +
  a_\delta \delta \wedge \alpha + a_\eta \eta \wedge \alpha$ and
  $a_\gamma \gamma \wedge \beta  +
  a_\delta \delta \wedge \beta + a_\eta \eta \wedge \beta$ are
contained in $[\ngo_1 , \ngo_1]$. Hence $a'_\alpha = a'_\beta =
0$, and so
  $v_1 = 0$.
  Thus if
 $v \in V_2$, then  $\ngo_1$ is contained in $V_2 \oplus V \oplus W$.
  Similarly we can prove that $\ngo_1$ is contained in $V \oplus V_1
  \oplus W$ or  $V \oplus V_2 \oplus W$.
  Hence $\ngo_1$ is contained in $V_i \oplus V_j \oplus W$ for some  $i , j
 \in \{1, 2 \}$. But then $[\ngo_1 , \ngo_1 ] = 0$ which is a
 contradiction. This proves the proposition.
\end{proof}

\end{document}